# ON THE TORIC ALGEBRA OF GRAPHICAL MODELS

By Dan Geiger, Christopher Meek and Bernd Sturmfels[1]

*Technion–Haifa, Microsoft Research and University of California, Berkeley*

We formulate necessary and sufficient conditions for an arbitrary discrete probability distribution to factor according to an undirected graphical model, or a log-linear model, or other more general exponential models. For decomposable graphical models these conditions are equivalent to a set of conditional independence statements similar to the Hammersley–Clifford theorem; however, we show that for nondecomposable graphical models they are not. We also show that nondecomposable models can have nonrational maximum likelihood estimates. These results are used to give several novel characterizations of decomposable graphical models.

**1. Introduction.** Exponential models for discrete data have a long history in statistics. In this paper, we take an algebraic approach to analyzing exponential models. Our starting point is to describe a class of exponential models for discrete distributions in terms of a polynomial mapping from a set of parameters to distributions. These models include two important well-known classes of models: the log-linear model and, an important type of log-linear model, the undirected graphical model. Representing the models as polynomials rather than in the more standard exponential representation allows us to use tools from computational algebraic geometry (e.g., [6]) to analyze the algebraic properties of these models.

We begin by providing necessary and sufficient conditions for a discrete probability distribution to factor according to an undirected graphical model, or a log-linear model, or a more general exponential model. The factorization of distributions according to these classes of models is well studied (see,

Received March 2003; revised July 2005.

[1]Supported in part by NSF Grant DMS-02-00729.

e.g., [2, 4, 5]). Unlike previous analyses that either assume positivity or use the exponential representation, our use of a polynomial representation allows us to provide a uniform treatment of the factorization for both positive and nonpositive distributions.

Next, utilizing computational tools and results from algebraic geometry, we analyze constraints imposed on distributions by specific exponential models. More specifically, we transform parametrically defined models into implicit descriptions. The following example illustrates the concept of an implicit description.

EXAMPLE 1.   Consider probability distributions over three binary variables $A, B, C$ defined parametrically as

$$P(a, b, c) = p_{abc} \propto \psi_{AB}(a, b)\psi_{BC}(b, c).$$

This can be viewed as either a log-linear model with generators $(AB)$ and $(BC)$ or as an undirected graphical model $A - B - C$. The corresponding implicit description is given as follows. A probability distribution $P = (p_{000}, p_{001}, p_{010}, p_{011}, p_{100}, p_{101}, p_{110}, p_{111})$ factors according to this model if and only if

$$p_{001}p_{100} = p_{000}p_{101} \quad \text{and} \quad p_{011}p_{110} = p_{010}p_{111}.$$

Perhaps the most well-known example of implicit descriptions of statistical models is given by the Hammersley–Clifford theorem (e.g., [3, 19]) which characterizes the factorization of strictly positive distributions with respect to undirected graphs.

Our analysis provides insight into two distinct but related approaches that have been used to study many types of graphical models including undirected and directed graphical models. The first approach is to define graphical models by specifying a graph according to which a probability distribution must factor in order to belong to the graphical model. This approach was emphasized, for example, by Darroch, Lauritzen and Speed [7]. The second approach is to define graphical models by specifying, through a graph, a set of conditional independence statements which a probability distribution must satisfy in order to belong to the graphical model. This direction was emphasized, for example, by Pearl [22] and Geiger and Pearl [13]. Lauritzen ([19], Chapter 3) compared these approaches and herein we extend his analysis. Our analysis allows us to identify the difference in these approaches and provides several novel characterizations of decomposable graphical models.

We note that using tools from computational algebra in the study of implicit descriptions of statistical models is not new. For instance, Settimi and Smith [26] and Geiger, Heckerman, King and Meek [12] analyze the geometric structure of directed graphical models with and without latent variables from the perspective of real algebraic geometry. In addition, Pistone, Riccomagno and Wynn [23] have used commutative algebra to study what we term the binary four-cycle model (see Example 4) and Garcia, Stillman and Sturmfels [11] have used commutative algebra to study graphical models with latent variables.

The paper is organized as follows. In Section 2 we define a class of exponential models and describe log-linear and undirected graphical models. In Section 3 we provide necessary and sufficient conditions for a discrete probability distribution to factor according to such an exponential model or to be the limit of distributions that factor. In Section 4 we focus our attention on undirected graphical models. We demonstrate that every nondecomposable model implies nonconditional-independence constraints, and show the possibility of nonrational maximum likelihood estimates for some nondecomposable models. Our analysis is summarized in Theorem 4.4, which provides various characterizations of decomposable models.

**2. Exponential, log-linear and graphical models.** Our objects of study are certain statistical models for a finite state space $\mathcal{X}$. We identify $\mathcal{X}$ with the set $\{1, 2, \ldots, m\}$ and define a probability distribution over $\mathcal{X}$ to be a vector $P = (p_1, \ldots, p_m)$ in $\mathbf{R}_{\geq 0}^m$ such that $p_1 + \cdots + p_m = 1$.

The class of models to be considered consists of discrete probability distributions defined via a $d \times m$ matrix $A = (a_{ij})$ of nonnegative integers. One technical assumption we will make about the matrix $A$ is that all its column sums are equal, that is, $\sum_{i=1}^{d} a_{i1} = \sum_{i=1}^{d} a_{i2} = \cdots = \sum_{i=1}^{d} a_{im}$. We say that a probability distribution $P$ *belongs to model* $A$ if and only if $P$ is in the image of the monomial mapping $\phi_A$ which takes nonnegative real $d$-vectors to nonnegative real $m$-vectors:

$$(2.1) \quad \phi_A : \mathbf{R}_{\geq 0}^d \to \mathbf{R}_{\geq 0}^m, \qquad (t_1, \ldots, t_d) \mapsto \left( \prod_i t_i^{a_{i1}}, \prod_i t_i^{a_{i2}}, \ldots, \prod_i t_i^{a_{im}} \right),$$

where, as we do throughout the paper, we adopt the convention that $t^0 = 1$ for $t \geq 0$. When $P$ belongs to model $A$, we also say that $P$ *factors* according to model $A$. The models described by (2.1) are usually described in the statistical literature as exponential families (models) of the form

$$(2.2) \quad P_\theta(x) = Z(\theta) e^{\langle \theta, T(x) \rangle}, \qquad \theta \in [-\infty, \infty)^d,$$

where $x \in \mathcal{X}$, $Z(\theta)$ is a normalizing constant, $\langle \cdot, \cdot \rangle$ denotes an inner product and *sufficient statistics* $T : \mathcal{X} \mapsto \mathbf{Z}^d \backslash \{\mathbf{0}\}$ where $\mathbf{Z}$ denotes the set of integers and $\mathbf{0}$ is a vector of $d$ zeroes.

The classes of models in (2.2) and (2.1) are identical. Note that each column of $A$ corresponds to a different state $x$ of $\mathcal{X}$. Thus, a model defined by (2.2) with sufficient statistics $T(x)$ is equivalent to a model defined by (2.1) with matrix $A$ if and only if the columns $a_j$ of $A$ coincide with the corresponding $T(x)$. For a particular distribution from a model of the form (2.2) with sufficient statistics $T(x)$ and parameters $\theta$, the corresponding parameters $t_i$ for the corresponding model in (2.1) are $t_i = \exp(\theta_i)$ where $\exp(-\infty) = 0$. We describe the models given by (2.2) in terms of a polynomial map because, as we shall see, this description allows us to use commutative algebra to provide algebraic descriptions of interesting properties of these models. Unlike the exponential representation, the use of a polynomial map allows us to directly analyze distributions that do not have full support. Recall that the *support* of an $m$-dimensional vector $v$ is the set of indices $\operatorname{supp}(v) = \{i \in \{1, \ldots, m\} : v_i \neq 0\}$.

This class of models includes log-linear and undirected graphical models used in the analysis of multiway contingency tables. When analyzing multiway contingency tables, the state space is a product space $\mathcal{X} = \prod_{X_j \in \mathbf{X}} I_{X_j}$ where $\mathbf{X} = \{X_1, \ldots, X_n\}$ is a set of (random) variables, called *factors*, and $I_{X_j}$ is the set of levels (or states) for the factor $X_j$. A log-linear model is defined by a collection $\mathcal{G} = \{\mathcal{G}_1, \ldots, \mathcal{G}_m\}$ of subsets of $\mathbf{X}$. We refer to the $\mathcal{G}_i$ as the *generators* of the log-linear model. A *log-linear model* for a set of generators $\mathcal{G}$ is defined as

$$P(x) \propto \prod_{\mathcal{G}_i \in \mathcal{G}} \psi_{\mathcal{G}_i}(x),$$

where $x \in \mathcal{X}$ is an instantiation of the variables in $\mathbf{X}$ and $\psi_{\mathcal{G}_i}(x)$ is a *potential function* that depends on $x$ only through the values of the variables in $\mathcal{G}_i$. This log-linear model can be represented in the following way by a matrix $A$ as in (2.1). The columns of $A$ are indexed by $\mathcal{X} = \prod_{X_j \in \mathbf{X}} I_{X_j}$. The rows of $A$ are indexed by pairs consisting of a generator $\mathcal{G}_i$ and an element of $\prod_{X_j \in \mathcal{G}_i} I_{X_j}$. All entries of $A$ are either zero or 1. The entry is 1 if and only if the element in the row index is equal to the projection of the column index to the factors in the generator of the row index.

EXAMPLE 2. The no-three-way interaction model for binary factors $X_1$, $X_2, X_3$ has generators $\mathcal{G} = \{\{X_1, X_2\}, \{X_2, X_3\}, \{X_1, X_3\}\}$ and is represented by the matrix

$$\begin{array}{c|cccccccc}
 & p_{000} & p_{001} & p_{010} & p_{011} & p_{100} & p_{101} & p_{110} & p_{111} \\
t_1 \equiv \psi_{\{1,2\}}(00) & 1 & 1 & 0 & 0 & 0 & 0 & 0 & 0 \\
t_2 \equiv \psi_{\{1,2\}}(01) & 0 & 0 & 1 & 1 & 0 & 0 & 0 & 0 \\
t_3 \equiv \psi_{\{1,2\}}(10) & 0 & 0 & 0 & 0 & 1 & 1 & 0 & 0 \\
t_4 \equiv \psi_{\{1,2\}}(11) & 0 & 0 & 0 & 0 & 0 & 0 & 1 & 1 \\
t_5 \equiv \psi_{\{2,3\}}(00) & 1 & 0 & 0 & 0 & 1 & 0 & 0 & 0 \\
t_6 \equiv \psi_{\{2,3\}}(01) & 0 & 1 & 0 & 0 & 0 & 1 & 0 & 0 \\
t_7 \equiv \psi_{\{2,3\}}(10) & 0 & 0 & 1 & 0 & 0 & 0 & 1 & 0 \\
t_8 \equiv \psi_{\{2,3\}}(11) & 0 & 0 & 0 & 1 & 0 & 0 & 0 & 1 \\
t_9 \equiv \psi_{\{1,3\}}(00) & 1 & 0 & 1 & 0 & 0 & 0 & 0 & 0 \\
t_{10} \equiv \psi_{\{1,3\}}(01) & 0 & 1 & 0 & 1 & 0 & 0 & 0 & 0 \\
t_{11} \equiv \psi_{\{1,3\}}(10) & 0 & 0 & 0 & 0 & 1 & 0 & 1 & 0 \\
t_{12} \equiv \psi_{\{1,3\}}(11) & 0 & 0 & 0 & 0 & 0 & 1 & 0 & 1
\end{array}.$$

A probability distribution $P = (p_{000}, p_{001}, p_{010}, p_{011}, p_{100}, p_{101}, p_{110}, p_{111})$ factors in the no-three-way interaction model if and only if it lies in the image of the associated monomial mapping

$$\begin{aligned}
(2.3) \quad & \phi_A : \mathbf{R}_{\geq 0}^{12} \to \mathbf{R}_{\geq 0}^8, \\
& (t_1, \dots, t_{12}) \to (t_1 t_5 t_9, t_1 t_6 t_{10}, t_2 t_7 t_9, t_2 t_8 t_{10}, t_3 t_5 t_{11}, t_3 t_6 t_{12}, t_4 t_7 t_{11}, t_4 t_8 t_{12}).
\end{aligned}$$

An important type of log-linear model is the undirected graphical model [19]. Such a model is specified by an undirected graph $G$ with vertex set $\mathbf{X}$ and edge set $\mathbf{E}$. The *undirected graphical model* for the graph $G$ is the log-linear model in which the generators are the *cliques* (maximal complete subgraphs) of the undirected graph $G$. The matrix $A$ of (2.1) is a function of the graph $G$ and we write it as $A(G)$. Example 2 shows a log-linear model that is not graphical.

EXAMPLE 3. The three-variable-chain undirected graphical model with graph $G$ equal to $X_1 - X_2 - X_3$ has generators $\mathcal{G} = \{\{X_1, X_2\}, \{X_2, X_3\}\}$. When each $X_i$ is a binary variable, the matrix $A(G)$ is identical to the first eight rows of the matrix of Example 2.

An undirected graphical model is said to be a *decomposable graphical model* if and only if the graph $G$ is chordal—that is, if every cycle of length 4 or more has a chord. The undirected graphical model given in Example 3 is a decomposable graphical model. We conclude this section with the four-cycle undirected graphical model, the simplest nondecomposable graphical model. This model will be examined in detail in Section 4.4.

EXAMPLE 4. The four-cycle undirected graphical model for binary variables with graph $G$ having four edges $X_1 - X_2$, $X_2 - X_3$, $X_3 - X_4$ and $X_1 - X_4$ has generators $\mathcal{G} = \{\{X_1, X_2\}, \{X_2, X_3\}, \{X_3, X_4\}, \{X_1, X_4\}\}$ and is represented by the following matrix $A(G)$:

| | $p_{0000}$ | $p_{0001}$ | $p_{0010}$ | $p_{0011}$ | $p_{0100}$ | $p_{0101}$ | $p_{0110}$ | $p_{0111}$ | $p_{1000}$ | $p_{1001}$ | $p_{1010}$ | $p_{1011}$ | $p_{1100}$ | $p_{1101}$ | $p_{1110}$ | $p_{1111}$ |
|---|---|---|---|---|---|---|---|---|---|---|---|---|---|---|---|---|
| $\psi_{\{1,2\}}(00)$ | 1 | 1 | 1 | 1 | 0 | 0 | 0 | 0 | 0 | 0 | 0 | 0 | 0 | 0 | 0 | 0 |
| $\psi_{\{1,2\}}(01)$ | 0 | 0 | 0 | 0 | 1 | 1 | 1 | 1 | 0 | 0 | 0 | 0 | 0 | 0 | 0 | 0 |
| $\psi_{\{1,2\}}(10)$ | 0 | 0 | 0 | 0 | 0 | 0 | 0 | 0 | 1 | 1 | 1 | 1 | 0 | 0 | 0 | 0 |
| $\psi_{\{1,2\}}(11)$ | 0 | 0 | 0 | 0 | 0 | 0 | 0 | 0 | 0 | 0 | 0 | 0 | 1 | 1 | 1 | 1 |
| $\psi_{\{2,3\}}(00)$ | 1 | 1 | 0 | 0 | 0 | 0 | 0 | 0 | 1 | 1 | 0 | 0 | 0 | 0 | 0 | 0 |
| $\psi_{\{2,3\}}(01)$ | 0 | 0 | 1 | 1 | 0 | 0 | 0 | 0 | 0 | 0 | 1 | 1 | 0 | 0 | 0 | 0 |
| $\psi_{\{2,3\}}(10)$ | 0 | 0 | 0 | 0 | 1 | 1 | 0 | 0 | 0 | 0 | 0 | 0 | 1 | 1 | 0 | 0 |
| $\psi_{\{2,3\}}(11)$ | 0 | 0 | 0 | 0 | 0 | 0 | 1 | 1 | 0 | 0 | 0 | 0 | 0 | 0 | 1 | 1 |
| $\psi_{\{3,4\}}(00)$ | 1 | 0 | 0 | 0 | 1 | 0 | 0 | 0 | 1 | 0 | 0 | 0 | 1 | 0 | 0 | 0 |
| $\psi_{\{3,4\}}(01)$ | 0 | 1 | 0 | 0 | 0 | 1 | 0 | 0 | 0 | 1 | 0 | 0 | 0 | 1 | 0 | 0 |
| $\psi_{\{3,4\}}(10)$ | 0 | 0 | 1 | 0 | 0 | 0 | 1 | 0 | 0 | 0 | 1 | 0 | 0 | 0 | 1 | 0 |
| $\psi_{\{3,4\}}(11)$ | 0 | 0 | 0 | 1 | 0 | 0 | 0 | 1 | 0 | 0 | 0 | 1 | 0 | 0 | 0 | 1 |
| $\psi_{\{1,4\}}(00)$ | 1 | 0 | 1 | 0 | 1 | 0 | 1 | 0 | 0 | 0 | 0 | 0 | 0 | 0 | 0 | 0 |
| $\psi_{\{1,4\}}(01)$ | 0 | 1 | 0 | 1 | 0 | 1 | 0 | 1 | 0 | 0 | 0 | 0 | 0 | 0 | 0 | 0 |
| $\psi_{\{1,4\}}(10)$ | 0 | 0 | 0 | 0 | 0 | 0 | 0 | 0 | 1 | 0 | 1 | 0 | 1 | 0 | 1 | 0 |
| $\psi_{\{1,4\}}(11)$ | 0 | 0 | 0 | 0 | 0 | 0 | 0 | 0 | 0 | 1 | 0 | 1 | 0 | 1 | 0 | 1 |

## 3. Exponential models and toric varieties.

In this section, we study the algebraic structure of the exponential models in (2.1). We provide a characterization of those distributions that factor according to a model and of those distributions that are the limit of distributions that factor. Finally, we describe how one can use tools from commutative algebra to obtain a complete description of the set of distributions that factor according to a model of the form (2.1) or are the limit of distributions that factor in terms of polynomial equations not involving model parameters.

3.1. *Distributions that factor and limits of distributions that factor.* We formulate necessary and sufficient conditions for a probability distribution to factor according to a matrix $A$ and for a distribution to be the limit of distributions that factor according to a matrix $A$.

The factorization of distributions according to exponential models is well studied (e.g., [2, 4, 5, 24]). Typically the analysis of these models is carried out using the exponential form given in (2.2). This type of analysis leads to the treatment of nonpositive distributions as special limiting cases such as the "boundaries at infinity" of Čencov [5]. The factorization of distributions according to log-linear models is also well studied (e.g., [8, 14, 15]). These analyses provide characterizations of factorization but only for positive distributions. By utilizing the product form representation of (2.1), we provide a uniform treatment of the factorization for both positive and nonpositive distributions.

Another alternative to our approach and to using an exponential representation is given in Lauritzen's [18] development of a generalization of exponential models called general exponential models. The general exponential model treats sufficient statistics as values in a commutative semigroup and replaces the exponential function with the members of a dual semigroup defined in terms of a homomorphism from the semigroup of sufficient statistics to the semigroup $(\mathbf{R}_{\geq 0}, \cdot)$. In certain examples, this approach yields a uniform treatment of positive and nonpositive distributions.

The characterization for distributions of the form of (2.1) is provided in terms of a condition on the support of the distribution and a set of algebraic constraints. We begin with the condition on the support of the distribution. Let $a_j = (a_{1j}, \ldots, a_{dj})$ denote the $j$th column vector of the $d \times m$ matrix $A$. Note that $\mathrm{supp}(a_j) \subseteq \{1, 2, \ldots, d\}$. A subset $F$ of $\{1, \ldots, m\}$ is said to be $A$-feasible if, for every $j \in \{1, \ldots, m\} \backslash F$, the support $\mathrm{supp}(a_j)$ of the vector $a_j$ is not contained in $\bigcup_{l \in F} \mathrm{supp}(a_l)$. Note that, trivially, the set $\{1, \ldots, m\}$ is $A$-feasible.

LEMMA 1. *A probability distribution $P$ factors according to $A$ only if the support of $P$ is $A$-feasible.*

PROOF. Let $P$ be a probability distribution which factors according to $A$, that is, $P \in \mathrm{image}(\phi_A)$. We must show that $F = \mathrm{supp}(P)$ is $A$-feasible. Let $(t_1, \ldots, t_d)$ be any preimage of $P$ under $\phi_A$. Then

$$(3.1) \qquad p_j = \begin{cases} \displaystyle\prod_{i=1}^{d} t_i^{a_{ij}} > 0, & \text{for } j \in F, \\ \displaystyle\prod_{i=1}^{d} t_i^{a_{ij}} = 0, & \text{for } j \notin F. \end{cases}$$

Suppose that $F$ is not $A$-feasible. Then $\mathrm{supp}(a_k)$ lies in $\bigcup_{l \in F} \mathrm{supp}(a_l)$ for some $k \notin F$. Consequently for every $i \in \mathrm{supp}(a_k)$, there exists an $f \in F$ such that $a_{if} > 0$. Hence, due to (3.1), $t_i > 0$ for every $i \in \mathrm{supp}(a_k)$. Thus $p_k = \prod_{i \in \mathrm{supp}(a_k)} t_i^{a_{ik}} > 0$ contrary to our assumption that $k \notin F$. □

Next we turn to the algebraic condition. The *nonnegative toric variety* $X_A$ is the set of all vectors $(x_1, \ldots, x_m) \in \mathbf{R}_{\geq 0}^{m}$ which satisfy

$$(3.2) \qquad x_1^{u_1} x_2^{u_2} \cdots x_m^{u_m} = x_1^{v_1} x_2^{v_2} \cdots x_m^{v_m},$$

whenever $u = (u_1, \ldots, u_m)$ and $v = (v_1, \ldots, v_m)$ are vectors of nonnegative integers which satisfy the $d$ linear relations

$$(3.3) \qquad u_1 a_1 + u_2 a_2 + \cdots + u_m a_m = v_1 a_1 + v_2 a_2 + \cdots + v_m a_m.$$

Note that (3.3) merely states that $u - v$ is in the kernel of the matrix $A$, that is, the matrix $A$ times the column vector $u - v$ is zero. Since the exponents $u_1, \ldots, u_m, v_1, \ldots, v_m$ used in (3.2) were assumed to be integers, the set $X_A$ is indeed an algebraic variety, that is, the zero set of a system of polynomial equations.

LEMMA 2. *A probability distribution $P$ factors according to $A$ only if $P$ lies in the nonnegative toric variety $X_A$.*

PROOF. We need to show that the image of $\phi_A$ is a subset of $X_A$. Indeed, suppose that $x = (x_1, \ldots, x_m) \in \text{image}(\phi_A)$. There exist nonnegative reals $t_1, \ldots, t_d$ such that $x_i = t_1^{a_{1i}} t_2^{a_{2i}} \cdots t_d^{a_{di}}$ for $i = 1, \ldots, m$. This implies that (3.2) has the form

$$\left( \prod_{j=1}^{d} t_j^{a_{j1}} \right)^{u_1} \left( \prod_{j=1}^{d} t_j^{a_{j2}} \right)^{u_2} \cdots \left( \prod_{j=1}^{d} t_j^{a_{jm}} \right)^{u_m}$$

$$= \left( \prod_{j=1}^{d} t_j^{a_{j1}} \right)^{v_1} \left( \prod_{j=1}^{d} t_j^{a_{j2}} \right)^{v_2} \cdots \left( \prod_{j=1}^{d} t_j^{a_{jm}} \right)^{v_m}$$

and hence it holds whenever (3.3) holds. Thus, $x$ lies in $X_A$. □

The following theorem provides a characterization of distributions that factor in terms of these two conditions.

THEOREM 3.1. *A probability distribution $P$ factors according to $A$ if and only if $P$ lies in the nonnegative toric variety $X_A$ and the support of $P$ is $A$-feasible.*

The only-if direction has been proved in Lemmas 1 and 2. The if direction is provided in the Appendix.

We now turn our discussion to the set of distributions that do not factor but are the limit of distributions that factor. In general, image($\phi_A$) is not a closed subset of the orthant $\mathbf{R}_{\geq 0}^m$. This is important because if there are distributions that do not factor according to a model but are the limit of distributions that do factor, then there is no unique maximum likelihood estimate (MLE) for some data sets. See Section 4 for an analysis of such phenomena in the four-cycle undirected graphical model. We will see in Theorem 4.4 that image($\phi_A$) is closed for an undirected graphical model if and only if the model is decomposable.

Our next theorem says that the set of probability distributions which lie in the toric variety $X_A$ coincides with those in the closure of the image of $\phi_A$—that is, $X_A = \text{closure}(\text{image}(\phi_A))$. Note that the closure can be taken in either the usual metric topology or in the Zariski topology because the closures of an image of a polynomial map taken in these topologies are the same. This result means that $P \in X_A$ if and only if $P$ factors according to $A$, or $P$ is the limit of probability distributions which factor according to $A$. The set of distributions in $X_A$, when $A$ consists only of zeroes and 1's, is called an extended log-linear model by Lauritzen [19]. Thus Theorem 3.2 below amounts to an algebraic description of extended exponential models and, thus, extended log-linear models and extended undirected graphical models.

THEOREM 3.2. *A probability distribution $P$ factors according to $A$ or is the limit of probability distributions that factor according to $A$ if and only if $P$ lies in the nonnegative toric variety $X_A$.*

The proof of Theorem 3.2 is provided in the Appendix.

Theorems 3.1 and 3.2 together characterize probability distributions in $X_A \setminus \text{image}(\phi_A)$, namely, distributions that do not factor but are the limit of distributions that do factor. These distributions are those that lie in $X_A$ but have a support which is not $A$-feasible.

3.2. *Describing exponential models by binomial equations.* In this section we describe an implicit representation of the toric variety that contains the distributions that factor according to the model in (2.1). The implicit representation is given in terms of the common zero set of a finite list of polynomial equations. These polynomial equations are interesting from both an algorithmic and a theoretical point of view in that they describe constraints on probability distributions that must hold for any distribution that factors according to the model.

Implicit representations of statistical models (see Example 1) are naturally described using the language of ideals and varieties. We briefly review these basic concepts from algebra and refer the reader to an excellent text by Cox, Little and O'Shea [6] for more details. All algebra terminology we use which is not defined in this paper can be found in [6].

We work in the ring $\mathbf{R}[x] = \mathbf{R}[x_1, \ldots, x_m]$ of polynomials with real coefficients in the indeterminates $x_1, \ldots, x_m$. An *ideal* $I$ is a nonempty subset of $\mathbf{R}[x]$ which satisfies two properties: (1) if $q_1, q_2 \in I$, then $q_1 + q_2 \in I$, and (2) if $b \in \mathbf{R}[x]$, and $q \in I$, then $bq \in I$. With every ideal $I$ in $\mathbf{R}[x]$ we associate a set of *varieties*,

$$X_I^K = \{x \in K^m : q(x) = 0 \text{ for every } q \in I\},$$

where $K$ denotes either the positive real numbers $\mathbf{R}_{>0}$ or the nonnegative real numbers $\mathbf{R}_{\geq 0}$. To simplify the notation, we write $X^{>0}$ and $X$ rather than $X^{\mathbf{R}_{>0}}$ and $X^{\mathbf{R}_{\geq 0}}$, respectively, and drop the explicit reference to the ideal when the associated ideal is clear from context. For $x \in K^m$, testing $x \in X^K$ is equivalent to checking that $q(x) = 0$ for all $q \in I$. In the analysis of statistical models $K = \mathbf{R}_{\geq 0}$ corresponds to the set of (nonnegative) probability distributions, and $K = \mathbf{R}_{>0}$ corresponds to the set of strictly positive probability distributions.

The task of checking that a point (e.g., a distribution) is in the zero set of each of the polynomials in an ideal appears extremely hard, but there are two fundamental results which make it more tractable. *Hilbert's basis theorem* states that every ideal in $\mathbf{R}[x]$ is *finally generated*, namely, every ideal $I$ in $\mathbf{R}[x]$ contains a finite subset $\{g_1, \ldots, g_n\}$, called an *ideal basis of $I$*, such that every $q \in I$ can be written as $q(x) = \sum_{i=1}^{n} b_i(x) g_i(x)$ where $b_i$ are polynomials in $\mathbf{R}[x]$. Consequently, a point $x$ in $K^m$ lies in $X^K$ if and only if $g_1(x) = \cdots = g_n(x) = 0$. The ideal generated by a set of polynomials $g = \{g_1, \ldots, g_n\}$ is denoted by $\langle g_1, \ldots, g_n \rangle$. The second fundamental result, by Buchberger, is an algorithm that produces a distinguished ideal basis, called a *Gröbner basis*, for any given ideal $I$. An ideal basis $g$ for $I$ is a Gröbner basis for $I$ in some term order (say lexicographical, or reverse lexicographical order) if the set of highest-ordered terms of the polynomials in $g$ generates the ideal generated by the highest-order terms of all polynomials in $I$.

The important property of Gröbner bases is that they allow one to check, in an efficient manner, whether a polynomial constraint belongs to an ideal. For example, if one obtains a small Gröbner basis for a graphical model under study, then one can use it to answer whether any cross product ratio, or any other polynomial constraint, must hold in that model. The focus on studying the ideals rather than the associated varieties also stems from the complexities introduced by allowing probability distributions that are not strictly positive.

In this paper, we consider ideals generated by a set of polynomials each having precisely two terms. Such polynomials are sometimes called *binomials*. The *toric ideal $I_A$* associated with a $d \times m$ integer matrix $A$ is generated by the binomials $x_1^{u_1} \cdots x_m^{u_m} - x_1^{v_1} \cdots x_m^{v_m}$ satisfying (3.3). A variety corresponding to a toric ideal is called a *toric* variety. An introduction to toric ideals can be found in [28]. We can now rewrite Theorems 3.1 and 3.2 as follows.

THEOREM 3.3. *A probability distribution $P$ factors according to an exponential model $A$ if and only if the support of $P$ is $A$-feasible and all polynomials in an ideal basis of the toric ideal $I_A$ vanish at $P$.*

THEOREM 3.4. *A probability distribution $P$ is the limit of probability distributions that factor according to $A$ if and only if all polynomials in an ideal basis of the toric ideal $I_A$ vanish at $P$.*

We call these the *factorization theorem* and the *limit factorization theorem*, respectively. Thus, if we know a small ideal basis for $I_A$, then we can efficiently test whether or not a distribution $P$ lies in $X_A$ by checking that $P$ satisfies these polynomials. It is important to note that it is frequently possible to replace the ideal basis $g$ of an ideal $I = \langle g_1, \ldots, g_n \rangle$ with a smaller basis $g'$ for an ideal $J$ such that the variety $X_I$ agrees with the variety $X_J$ in the nonnegative orthant (i.e., $X_I^{\mathbf{R}_{\geq 0}} = X_J^{\mathbf{R}_{\geq 0}}$). Thus, one can often identify smaller sets than the ideal basis for $I_A$ for use in Theorems 3.3 and 3.4 when testing a distribution. We demonstrate that, even in the case of decomposable models, the ideal basis for $I_A$ is typically larger than an ideal basis for an ideal whose zero set defines $X_A$. For arbitrary undirected graphical models, the Hammersley–Clifford theorem, to be discussed in the next section, defines a small subset of binomials whose zero set defines $X_{A(G)}^{>0}$, while it is an open problem to describe the Gröbner basis for an arbitrary undirected graphical model $G$.

**4. Algebraic analysis of graphical models.** The algebraic tools developed in Section 3 will now be applied to the undirected graphical models $A(G)$. We compare and contrast the Hammersley–Clifford theorem (e.g., [19], page 36; [3]) with the factorization Theorem 3.3 and limit factorization Theorem 3.4. We investigate the form of the ideal bases for decomposable and noncomposable models. We also study the algebraic complexity of the maximum likelihood estimator for undirected graphical models. Our main result is a characterization of decomposable graphical models in terms of their ideal basis, the rationality of maximum likelihood estimates, and whether the model contains all of its limit points.

4.1. *Quadratic polynomials representing conditional independence.* The set of probability distributions that satisfy a conditional independence statement can be regarded as an algebraic variety. In this subsection we explain how to derive the defining ideal of such a variety. The ideal basis will consist of certain quadratic polynomials which we call *cross-product differences* (CPDs). Given three discrete random variables $X, Y, Z$, we define

$$\begin{aligned}
(4.1) \quad & \mathrm{cpd}(X = \{x, x'\}, Y = \{y, y'\} | Z = z) \\
& \doteq P(x, y, z)P(x', y', z) - P(x', y, z)P(x, y', z),
\end{aligned}$$

where $x$ and $x'$ are levels of $X$ and $y$ and $y'$ are levels of $Y$ and $z$ is a level of $Z$. Note that cross-product differences are closely related to *cross-product ratios* (CPRs); the CPR is defined as follows:

$$(4.2) \quad \mathrm{cpr}(X = \{x, x'\}, Y = \{y, y'\} | Z = z) \doteq \frac{P(x, y, z)P(x', y', z)}{P(x', y, z)P(x, y', z)}.$$

Cross-product ratios (also called conditional odds-ratios) are a fundamental measure of association and interaction and are often used to interpret the parameters of a log-linear model (see, e.g., [1]). A CPD and the corresponding CPR constraint are identical in the sense that

$$\mathrm{cpd}(X = \{x, x'\}, Y = \{y, y'\} | Z = z) = 0 \qquad \text{if and only if}$$

$$\mathrm{cpr}(X = \{x, x'\}, Y = \{y, y'\} | Z = z) = 1,$$

provided the denominators in (4.2) are nonzero. We prefer CPD constraints to avoid dividing by zero for nonpositive distributions. However, when interpreting higher-degree binomials in the toric ideal of an undirected graphical model, it is convenient to describe the constraints in terms of cross-product ratios which can then be converted into binomial constraints by clearing the denominator as above. To simplify notation, when $X$ and $Y$ each represent a single binary variable, we write

$$(4.3) \qquad \mathrm{cpr}(X, Y | Z = z) = \frac{P(x, y, z) P(x', y', z)}{P(x', y, z) P(x, y', z)}.$$

Let $X_1, \ldots, X_n$ denote discrete variables, where $I_{X_j}$ is the set of levels of the variable $X_j$. We fix the polynomial ring $\mathbf{R}[\mathcal{X}]$ whose indeterminates are elementary probabilities $p_{a_1 a_2 \cdots a_n}$ which are indexed by the elements of $\mathcal{X} = I_{X_1} \times I_{X_2} \times \cdots \times I_{X_n}$. Conditional independence statements have the form

$$(4.4) \qquad\qquad X \text{ is independent of } Y \text{ given } Z,$$

where $X$, $Y$ and $Z$ are pairwise disjoint subsets of $\{X_1, \ldots, X_n\}$. The statement (4.4) translates into a large set of CPDs of the form (4.1). Namely, we take $\mathrm{cpd}(X = \{x, x'\}, Y = \{y, y'\} | Z = z)$, where $x, x'$ runs over distinct elements in $\prod_{X_i \in X} I_{X_i}$, where $y, y'$ runs over distinct elements in $\prod_{X_j \in Y} I_{X_j}$ and where $z$ runs over $\prod_{X_k \in Z} I_{X_k}$. Note that some of these CPDs may be redundant.

Each probability $P(x, y, z)$ occurring in the CPDs of a conditional independence statement is obtained by marginalizing over all of the elementary probabilities $p_{a_1 a_2 \cdots a_n}$ for which the indices agree with $x, y$ and $z$. This means that the probability $P(x, y, z)$ is a polynomial of degree 1 in $\mathbf{R}[\mathcal{X}]$. The linearity of probabilities and the form of the CPD in (4.1) lead to the following remark.

REMARK 1. Conditional independence statements translate into a system of CPDs that correspond to quadratic polynomials.

The conditional independence statement (4.4) is said to be *saturated* if $X \cup Y \cup Z = \{X_1, \ldots, X_n\}$. The fact that probabilities in the CPDs associated with a saturated conditional independence statement do not require marginalization leads to the following remark.

REMARK 2. Saturated conditional independence statements about the variables $X_1, \ldots, X_n$ translate into a set of quadratic binomials.

4.2. *Undirected graphical models and the Hammersley–Clifford theorem.* In this section, we describe sets of conditional independence statements derived from separation statements in undirected graphs and their connection to the factorization of distributions. Of particular interest is the Hammersley–Clifford theorem that relates the factorization of a strictly positive distribution $P$ according to an undirected graphical model to a set of conditional independence statements that must hold in $P$. We describe the Hammersley–Clifford theorem in the language of ideals and varieties and compare it to our factorization theorem.

Let $G$ be an undirected graphical model with variables $\{X_1, \ldots, X_n\}$ as before. We define $I_{\mathrm{pairwise}(G)}$ to be the ideal in $\mathbf{R}[\mathcal{X}]$ generated by the quadratic binomials corresponding to all the saturated conditional independence statements

(4.5)     $X_i$ is independent of $X_j$ given $\{X_1, \ldots, X_n\} \backslash \{X_i, X_j\}$,

where $(X_i, X_j)$ runs over all nonedges of the graph $G$. Note that (4.5) is saturated, so the polynomials arising from the construction in the previous section are indeed binomials. The ideal $I_{\mathrm{pairwise}(G)}$ defines a variety $X^K_{\mathrm{pairwise}(G)}$ where $K$ can be either $\mathbf{R}_{\geq 0}$ or $\mathbf{R}_{>0}$. When $K = \mathbf{R}_{\geq 0}$, the superscript of $X$ is dropped.

The *pairwise Markov property* is discussed in Section 3.2.1 of [19]. Lauritzen uses the notation $M_P(\mathcal{G})$ to denote the variety $X_{\mathrm{pairwise}(G)}$. We will also need the (saturated) *global Markov property*. This is described in our language as follows. We define $I_{\mathrm{global}(G)}$ to be the ideal in $\mathbf{R}[\mathcal{X}]$ generated by the quadratic binomials corresponding to all the saturated conditional independence statements (4.4) where *Z separates X from Y* in the graph $G$. (The term global Markov is often used to describe the set of conditional independence statements that follow from all separation statements rather than only the saturated separation statements. The fact that only the saturated statements are needed follows from simple properties of undirected graphs and conditional independence. The required conditional independence properties are properties C1 and C2 of [19], page 29. The required graph property is that any unsaturated separation statement in a graph is implied by a saturated separation fact also true in the graph.) This separation condition means that every path from a vertex in $X$ to a vertex in $Y$ must pass through some vertex in $Z$. The ideal $I_{\mathrm{global}(G)}$ defines a variety $X^K_{\mathrm{pairwise}(G)}$ where $K$ is either $\mathbf{R}_{\geq 0}$ or $\mathbf{R}_{>0}$. When $K = \mathbf{R}_{\geq 0}$, the superscript of $X$ is dropped. Lauritzen [19] states the following three inclusions, which hold for every graph $G$. Each of the following three inclusions can be strict:

(4.6)     $\mathrm{image}(\phi_{A(G)}) \subseteq X_{A(G)} \subseteq X_{\mathrm{global}(G)} \subseteq X_{\mathrm{pairwise}(G)}.$

The following example provides an illustration of quadratic polynomials generating $I_{\mathrm{pairwise}(G)}$.

EXAMPLE 5. Consider the four-cycle undirected graphical model of Example 4. This graph has four maximal cliques, one for each edge. The probability distributions defined by this model have the form

$$P(x_1, x_2, x_3, x_4) \propto \psi_{\{1,2\}}(x_1, x_2)\psi_{\{2,3\}}(x_2, x_3)\psi_{\{3,4\}}(x_3, x_4)\psi_{\{1,4\}}(x_1, x_4).$$
(4.7)

If all four variables are binary, then the pairwise ideal is

(4.8)
$$\begin{aligned}
I_{\mathrm{pairwise}(G)} = \langle &p_{1011}p_{1110} - p_{1010}p_{1111}, p_{0111}p_{1101} - p_{0101}p_{1111}, \\
&p_{1001}p_{1100} - p_{1000}p_{1101}, p_{0110}p_{1100} - p_{0100}p_{1110}, \\
&p_{0011}p_{1001} - p_{0001}p_{1011}, p_{0011}p_{0110} - p_{0010}p_{0111}, \\
&p_{0001}p_{0100} - p_{0000}p_{0101}, p_{0010}p_{1000} - p_{0000}p_{1010}\rangle.
\end{aligned}$$

This is a binomial ideal in a polynomial ring in sixteen indeterminates:

$$I_{\mathrm{pairwise}(G)} \subset \mathbf{R}[\mathcal{X}] = \mathbf{R}[p_{0000}, p_{0001}, p_{0010}, \ldots, p_{1111}].$$

The left column of four binomials in (4.8) represents the statement "$X_2$ is independent of $X_4$ given $\{X_1, X_3\}$," and the right column of four binomials in (4.8) represents the statement "$X_1$ is independent of $X_3$ given $\{X_2, X_4\}$." The variety $X_{\mathrm{pairwise}(G)}$ is the set of all points in $K^{16}$ which are common zeros of these eight binomials. Note that $I_{\mathrm{pairwise}(G)} = I_{\mathrm{global}(G)}$ for the four-cycle model and therefore, for this model, the right inclusion of (4.6) is an equality.

The following well-known theorem (e.g., [19], page 36) relates the ideal of pairwise conditional independence statements and factorization.

THEOREM 4.1 (Hammersley–Clifford). *Let $G$ be an undirected graphical model. A strictly positive probability distribution $P$ factors according to $A(G)$ if and only if $P$ is in the variety $X_{\mathrm{pairwise}(G)}^{>0}$; that is, $X_{A(G)}^{>0} = X_{\mathrm{pairwise}(G)}^{>0}$.*

Our factorization Theorem 3.3 generalizes the Hammersley–Clifford theorem in two respects. First, it does not require the probability distribution $P$ to be strictly positive. Second, it does not require the model represented by matrix $A$ to be an undirected graphical model. The main advantage of the Hammersley–Clifford theorem over the factorization theorem is computational. That is, the set $I_{\mathrm{pairwise}(G)}$ is easily described in terms of the structure of the graph while one must usually resort to a symbolic algebra program to produce an ideal basis or a Gröbner basis for $I_A$.

The proof of the Hammersley–Clifford theorem given in [19] actually establishes the following slightly stronger result: any integer vector in the kernel of the matrix $A(G)$ is an integer linear combination of the vectors $u - v$ corresponding to the binomials $p^u - p^v$ arising from the conditional independence statements for the nonadjacent pairs $(X_i, X_j)$ in $G$. Translating this statement from the additive notation into multiplicative notation, we obtain the following:

> A binomial $p^u - p^v$ lies in the toric ideal $I_{A(G)}$ of an undirected graphical model $A(G)$ if and only if some monomial multiple of it, that is, a binomial of the form $p^{u+w} - p^{v+w}$, lies in $I_{\mathrm{pairwise}(G)}$.

This fact is important for computational purposes. It means that we can use the quadratic binomials in $I_{\mathrm{pairwise}(G)}$ as input when computing the toric ideal $I_{A(G)}$ by Algorithm 12.3 of [28].

### 4.3. Decomposable models.
In this section, we discuss factorization and ideal bases for the variety of probability distributions corresponding to decomposable graphical models.

THEOREM 4.2 ([19], Proposition 3.19). *Let $G$ be a decomposable graphical model. A probability distribution $P$ factors according to $A(G)$ if and only if $P$ is in $X_{\mathrm{global}(G)}$.*

This theorem is analogous to the Hammersley–Clifford theorem in that it provides an implicit description of distributions that factor according to a decomposable graph in terms of conditional independence statements. Unlike the Hammersley–Clifford theorem, this theorem is not restricted to positive distributions. An immediate corollary to this theorem is the following: if $P$ is a limit of probability distributions that factor according to $A(G)$, then $P$ itself factors according to $A(G)$. This implies that the support of a distribution $P$ need not be tested in order to decide whether $P$ factors according to a decomposable model. Furthermore, for a decomposable graphical model $G$, two of the inclusions in (4.6) are equalities,

$$(4.9) \qquad \mathrm{image}(\phi_{A(G)}) = X_{A(G)} = X_{\mathrm{global}(G)} \subseteq X_{\mathrm{pairwise}(G)},$$

but the inclusion on the right-hand side is generally strict. The two equalities on the left are equivalent to Theorem 4.2. We shall see in Example 6 below that the inclusion on the right is strict for the four-chain model.

Not every toric ideal $I_A$ which is generated by quadratic binomials has a Gröbner basis consisting of quadratic binomials (see, e.g., [28]). It turns out that toric ideals arising from decomposable graphical models are well behaved in this regard.

THEOREM 4.3.    *Let $G$ be a decomposable graphical model. Then the set of quadratic binomials representing CPDs for saturated conditional independence statements for $G$ forms a Gröbner basis of the toric ideal $I_{A(G)}$.*

A nice proof of this theorem was given by Hoşten and Sullivant [17]. Their Theorem 4.17 explicitly constructs a minimal (and reduced) Gröbner basis for an arbitrary decomposable graphical model. From an algebraic point of view, we note that equation $X_{A(G)}^K = X_{\text{global}(G)}^K$ for chordal graphs $G$ holds not just for $K = \mathbf{R}_{\geq 0}$, but also for $K = \mathbf{R}$ and for $K = \mathbf{C}$. Takken [29] and Dobra [10] proved that this equality holds in the ideal-theoretic sense, namely, that $I_{A(G)} = I_{\text{global}(G)}$. This means that the CPDs (quadratic binomials) representing global conditional independence statements contain an ideal basis for the toric ideal $I_{A(G)}$ where $G$ is decomposable.

Some of the statistical implications of this result are explicated in Diaconis and Sturmfels [9]. They showed that every minimal ideal basis of the toric ideal $I_A$ provides a set of moves for a Markov chain Monte Carlo approach to sampling from the conditional distribution of data given sufficient statistics for discrete exponential families of the form (2.1). They showed that a minimal Gröbner basis guarantees that the resulting Markov chain is connected, and that no proper subset of such an ideal basis has this property.

We complete this section with an example that illustrates that the rightmost subset relation in (4.9) is strict and the fact that the Hammersley–Clifford theorem fails for nonpositive graphical models even when the models are decomposable.

EXAMPLE 6.    Consider the chain model $G_4$ for four binary variables $X_1 - X_2 - X_3 - X_4$. The ideal representing the pairwise Markov property is generated by twelve quadratic binomials:

$$
\begin{aligned}
I_{\text{pairwise}(G_4)} = \langle\, & p_{0010}p_{1000} - p_{0000}p_{1010},\, p_{0001}p_{1000} - p_{0000}p_{1001}, \\
& p_{0001}p_{0100} - p_{0000}p_{0101},\, p_{0011}p_{1001} - p_{0001}p_{1011}, \\
& p_{0011}p_{1010} - p_{0010}p_{1011},\, p_{0011}p_{0110} - p_{0010}p_{0111}, \\
& p_{0110}p_{1100} - p_{0100}p_{1110},\, p_{0101}p_{1100} - p_{0100}p_{1101}, \\
& p_{1001}p_{1100} - p_{1000}p_{1101},\, p_{0111}p_{1101} - p_{0101}p_{1111}, \\
& p_{0111}p_{1110} - p_{0110}p_{1111},\, p_{1011}p_{1110} - p_{1010}p_{1111}\,\rangle.
\end{aligned}
$$

There are many probability distributions which show that the inclusion in (4.9) is strict for this example. For instance, take $p_{0010} = p_{1111} = 1/2$ and all other 14 indeterminates zero. The twelve ideal generators of $I_{\text{pairwise}(G_4)}$ all vanish at this distribution but the binomial $p_{0011}p_{1110} - p_{0010}p_{1111} \in I_{A(G_4)}$ that is implied by the independence of $X_4$ and $\{X_1, X_2\}$ given $X_3$ does not.

4.4. *Nondecomposable models.* We now discuss nondecomposable undirected graphical models from the perspective of the factorization Theorem 3.3 and study the implicit description of distributions that factor in terms of the ideal bases for the toric ideal $I_{A(G)}$. These ideal bases contain polynomials which do not correspond to conditional independence statements. We explicitly describe the nonconditional-independence polynomials for the four-cycle model and demonstrate that the degree of the polynomial constraints describing factorization can grow exponentially in the number of variables.

Probability distributions which factor according to the four-cycle model of Example 4 must satisfy not just the eight quadratic binomials in (4.8), which arise from pairwise conditional independence statements, but they must satisfy certain additional polynomials of degree 4 listed in (4.10).

PROPOSITION 1. *Consider the four-cycle undirected graphical model of Example 4 with graph $G'$. A probability distribution $P$ factors according to the four-cycle or is the limit of probability distributions that factor according to the four-cycle if and only if $P$ satisfies the following ideal basis of the toric ideal $I_{A(G')}$:*

$$I_{A(G')} = I_{\mathrm{pairwise}(G')} + \langle f_{12}^{\mathrm{diff}}, f_{23}^{\mathrm{diff}}, f_{34}^{\mathrm{diff}}, f_{14}^{\mathrm{diff}}, f_{12}^{\mathrm{same}}, f_{23}^{\mathrm{same}}, f_{34}^{\mathrm{same}}, f_{14}^{\mathrm{same}} \rangle,$$

*where*

$$
\begin{aligned}
f_{12}^{\mathrm{diff}} &= p_{0100}p_{0111}p_{1001}p_{1010} - p_{0101}p_{0110}p_{1000}p_{1011}, \\
f_{23}^{\mathrm{diff}} &= p_{0010}p_{0101}p_{1011}p_{1100} - p_{0011}p_{0100}p_{1010}p_{1101}, \\
f_{34}^{\mathrm{diff}} &= p_{0001}p_{0110}p_{1010}p_{1101} - p_{0010}p_{0101}p_{1001}p_{1110}, \\
f_{14}^{\mathrm{diff}} &= p_{0001}p_{0111}p_{1010}p_{1100} - p_{0011}p_{0101}p_{1000}p_{1110}, \\
(4.10) \quad f_{12}^{\mathrm{same}} &= p_{0000}p_{0011}p_{1101}p_{1110} - p_{0001}p_{0010}p_{1100}p_{1111}, \\
f_{23}^{\mathrm{same}} &= p_{0000}p_{0111}p_{1001}p_{1110} - p_{0001}p_{0110}p_{1000}p_{1111}, \\
f_{34}^{\mathrm{same}} &= p_{0000}p_{0111}p_{1011}p_{1100} - p_{0011}p_{0100}p_{1000}p_{1111}, \\
f_{14}^{\mathrm{same}} &= p_{0000}p_{0110}p_{1011}p_{1101} - p_{0010}p_{0100}p_{1001}p_{1111}.
\end{aligned}
$$

The basis given in this proposition is obtained from Algorithm 12.3 of [28] using the eight quadratic generators of $I_{\mathrm{pairwise}(G')}$ and the polynomial map $\phi_{A(G')}$.

Next we provide an interpretation of the ideal basis of the four-cycle given in Proposition 1. The basis prescribed by (4.10) can be described in terms of a ratio of cross-product ratios. In particular, using the definition of CPR in (4.3), the eight new basis elements (4.10) can be written as follows:

$$\operatorname{cpr}(X_3, X_4 | X_1 X_2 = 01) / \operatorname{cpr}(X_3, X_4 | X_1 X_2 = 10) = 1,$$

$$\operatorname{cpr}(X_1, X_4 | X_2 X_3 = 01) / \operatorname{cpr}(X_1, X_4 | X_2 X_3 = 10) = 1,$$

$$\operatorname{cpr}(X_1, X_2 | X_3 X_4 = 01) / \operatorname{cpr}(X_1, X_2 | X_3 X_4 = 10) = 1,$$

$$\operatorname{cpr}(X_2, X_3 | X_1 X_4 = 01) / \operatorname{cpr}(X_2, X_3 | X_1 X_4 = 10) = 1,$$

(4.11)
$$\operatorname{cpr}(X_3, X_4 | X_1 X_2 = 00) / \operatorname{cpr}(X_3, X_4 | X_1 X_2 = 11) = 1,$$

$$\operatorname{cpr}(X_1, X_4 | X_2 X_3 = 00) / \operatorname{cpr}(X_1, X_4 | X_2 X_3 = 11) = 1,$$

$$\operatorname{cpr}(X_1, X_2 | X_3 X_4 = 00) / \operatorname{cpr}(X_1, X_2 | X_3 X_4 = 11) = 1,$$

$$\operatorname{cpr}(X_2, X_3 | X_1 X_4 = 00) / \operatorname{cpr}(X_2, X_3 | X_1 X_4 = 11) = 1.$$

These constraints force the association between adjacent variables in the four-cycle to be identical for various values of the remaining variables. Thus, these constraints, when written as polynomials rather than ratios of polynomials, are restricting higher-order interactions, but, surprisingly, are only needed for the characterization of nonpositive distributions.

Proposition 1 provides an ideal basis for the four-cycle undirected graphical model; however, the problem of explicitly providing a basis for an arbitrary undirected graphical model remains open.

We note that there is no general upper bound for the degrees of the binomials in the ideal basis of an undirected graphical model. For instance, if each variable in the four-cycle model has $p$ levels, then there exists a minimal generator of degree $\geq p$. Such a binomial can be derived from Proposition 14.14 in [28]. The next proposition demonstrates that the maximal degree of the polynomials in the ideal basis is unbounded when the complexity of the model increases even when all variables remain binary.

PROPOSITION 2. *There exists an undirected graphical model for $2n$ binary variables $X_1, \ldots, X_{2n}$ whose ideal basis contains a binomial of degree $2^n$.*

PROOF. Let $G$ be the undirected graphical model whose only nonedges are $\{X_i, X_{i+n}\}$ for $i = 1, 2, \ldots, n$. Thus this model represents $n$ pairs of noninteracting binary variables. Let $p^u$ denote the product of all indeterminates $p_{i_1 \cdots i_{2n}}$ such that $i_1 = i_3 = i_5 = \cdots = i_{2n-1}$ and $i_1$ has the same parity as $i_2 + i_4 + i_6 + \cdots + i_{2n}$, and let $p^v$ denote the product of all indeterminates $p_{i_1 \cdots i_{2n}}$ such that $i_1 = i_3 = i_5 = \cdots = i_{2n-1}$ and $i_1$ has parity different from $i_2 + i_4 + i_6 + \cdots + i_{2n}$. Then $p^u - p^v$ is a binomial of degree $2^n$ which lies in the toric ideal $I_{A(G)}$. It can be checked, for instance using Corollary 12.13 in [28], that $p^u - p^v$ is a minimal generator of $I_{A(G)}$. □

The undirected graphical models in the previous proof provide an interesting family for further study. Note that for $n = 2$ this is precisely the four-cycle model, and for $n = 3$ this is the edge graph of the octahedron, with cliques $\{1, 2, 3\}, \{1, 2, 6\}, \{1, 3, 5\}, \{1, 5, 6\}, \{2, 3, 4\}, \{2, 4, 6\}, \{3, 4, 5\}, \{4, 5, 6\}$. Here the binomial constructed in the proof of Proposition 2 equals

$$p^u - p^v = p_{000000}p_{000101}p_{010010}p_{010100}p_{101011}p_{101110}p_{111010}p_{111111}$$

$$- p_{000001}p_{000100}p_{010000}p_{010101}p_{101010}p_{101111}p_{111011}p_{111110}$$

which can also be written as a ratio of ratios of CPRs.

4.5. *Variety differences for the four-cycle model.* We focus on the following relationships which hold for every undirected graphical model [see (4.6)]:

$$\text{image}(\phi_{A(G)}) \subseteq X_{A(G)} = \text{closure}(\text{image}(\phi_{A(G)})) \subseteq X_{\text{global}(G)}.$$

Lauritzen showed via examples that both inclusions are strict for the four-cycle model, in contrast to (4.9) for decomposable models. In this section, we continue our algebraic analysis of the four-cycle model, studying the set differences $X_{A(G)} \setminus \text{image}(\phi_{A(G)})$ and $X_{\text{global}(G)} \setminus X_{A(G)}$. The examples considered herein are used in the proof of our characterization theorem of decomposable models (Theorem 4.4).

The distributions that lie in $X_{A(G)} \setminus \text{image}(\phi_{A(G)})$ are those that have a support which is not $A$-feasible. The following example from [19], page 37, illustrates such a distribution and is due to Moussouris [21].

EXAMPLE 7. Consider the probability distribution over four binary variables $X_1, X_2, X_3, X_4$ where

$$(4.12) \quad p_{0000} = p_{0001} = p_{1000} = p_{0011} = p_{1100} = p_{0111} = p_{1110} = p_{1111} = 1/8.$$

This distribution satisfies all 16 binomial generators of $I_{A(G')}$ where $A(G')$ is the $16 \times 16$ matrix in Example 4, and hence lies in the toric variety $X_{A(G')}$. However, this distribution does not factor according to the four-cycle because the support is not $A$-feasible. This can be seen from the matrix $A(G')$: if $F$ is the set of eight column indices appearing in (4.12), then $\bigcup_{l \in F} \text{supp}(a_l)$ consists of all 16 row indices of $A(G')$.

Because the distribution (4.12) is in $X_{A(G')}$ we know that it is the limit of distributions that factor. Lauritzen proves this by writing it explicitly as a limit of distributions that factor according to $G'$. This example highlights the importance of being $A$-feasible when it comes to factorization, and it illustrates our characterization of distributions that do not factor but are the limit of distributions that do factor, which is provided by Theorems 3.1 and 3.2.

We now discuss the set difference $X_{\text{global}(G)} \setminus X_{A(G)}$ of distributions that satisfy the global Markov property but are not limits of distributions that factor.

EXAMPLE 8.  Let $P$ be the distribution over four binary variables in which

$$p_{0100} = p_{0111} = p_{1001} = p_{1010} = 1/4.$$

This distribution satisfies the global Markov property for the four-cycle [i.e., it lies in $X_{\text{global}(G')}$]. However, it is not the limit of distributions that factor [i.e., it does not lie in $X_{A(G')}$ because $f_{12}^{\text{diff}}(P) \neq 0$]. Note that the other 15 generators of $X_{A(G')}$ vanish at $P$.

We note that probability distributions with these properties in the four-cycle models (for ternary variables) were found by Matúš and Studený [20]. Also see Example 3.15 in [19], page 41. This example also demonstrates an immediate corollary of Proposition 1 and the limit factorization Theorem 3.4: the set $X_{\text{pairwise}(G')} \setminus X_{A(G')}$ contains all probability distributions that do not factor and are not the limit of distributions that factor but satisfy the pairwise conditional independence statements of the Hammersley–Clifford theorem.

Finally, our contribution to the study of this four-cycle model is to provide a completely general algebraic method for describing the set $X_{\text{global}(G')} \setminus X_{A(G')}$. This set consists of all distributions in $X_{\text{global}(G')}$ except those which violate at least one polynomial in Proposition 1. For this claim to hold we need to show that none of the 16 generators listed in Proposition 1 for the four-cycle model is redundant in the limit factorization Theorem 3.4, that is, for any of these 16 binomials in the ideal basis of $I_{A(G)}$ there exists a probability distribution which satisfies the other 15 binomials but does not lie in $X_{A(G)}$. Example 8 provides one such distribution and others can be constructed in an analogous fashion.

4.6. *Maximum likelihood estimation.*  In this section, we consider the problem of maximum likelihood estimation for undirected graphical models. One of the nice properties of decomposable models is that the maximum likelihood estimates are provided by a simple ratio of counts (see, e.g., [19], page 91). We demonstrate that the situation with nondecomposable models and the general exponential models of Section 3 is not so nice.

For nondecomposable models the maximum likelihood estimate need not exist for the model as defined by (2.1). This can be seen by considering the problem of finding the maximum likelihood estimate (MLE) for the four-cycle undirected graphical model with four binary variables when given a data set with the empirical distribution given in Example 7. As we have seen, the support of this distribution is not $A$-feasible and, thus, cannot be parameterized as described in Example 4. Given that the distribution satisfies the generators of the ideal basis, we know that the distribution lies in the closure of the set of distributions parameterized in Example 4. This example demonstrates that MLE can fail to exist. In particular, if the empirical distribution is in the boundary of the model but cannot be factored according to the model, then the MLE will fail to exist.

It is natural to extend the model, as described in Section 3, to include the distributions that are limits of distributions that factor according to the model. For the remainder of the paper, we consider only extended undirected graphical models and extended log-linear models. This approach was used by Lauritzen ([19], Chapter 4) to demonstrate that the MLE always exists for extended log-linear models. As noted in Section 3, the toric variety for the model and its ideal provide algebraic descriptions of extended log-linear models. In fact, one can compute the MLE by using a purely algebraic approach by (1) parameterizing the model with cell counts, (2) forcing the set of polynomial generators in the ideal basis to be equal to zero, and (3) forcing the set of marginal counts for each of the possible values for the cliques of the undirected graphical model (or, more generally, the generators of the log-linear model) to match the sum of the associated cell counts. The MLE is the unique real-valued nonnegative solution to this set of polynomial equations (see, e.g., [19]).

Framing the problem of identifying the MLE as an algebraic problem allows the use of algebraic tools to analyze properties of the MLE for nondecomposable models. In the remainder of this section we use algebraic methods from Galois theory to demonstrate that the MLE for a nondecomposable model is not necessarily rational and that one cannot generally write the MLE for nondecomposable models in closed form.

Consider the four-cycle model for four binary variables (Example 4). We present the maximum likelihood estimation for this model in full detail for one explicit nontrivial data set, namely,

$$(4.13) \quad \begin{pmatrix} m_{0000} & m_{0001} & m_{0010} & m_{0011} \\ m_{0100} & m_{0101} & m_{0110} & m_{0111} \\ m_{1000} & m_{1001} & m_{1010} & m_{1011} \\ m_{1100} & m_{1101} & m_{1110} & m_{1111} \end{pmatrix} = \begin{pmatrix} 1 & 1 & 1 & 1 \\ 1 & 1 & 1 & 1 \\ 1 & 1 & 1 & 2 \\ 0 & 0 & 0 & 0 \end{pmatrix},$$

where $m_{ijkl}$ is the count of cases in which $X_1 = i$, $X_2 = j$, $X_3 = k$ and $X_4 = l$. The maximum likelihood estimate for our data set is a solution to a system of algebraic equations in 16 indeterminates $\hat{m}_{ijkl}$. The last four coordinates of the maximum likelihood estimate will automatically be zero,

$$\hat{m}_{1100} = \hat{m}_{1101} = \hat{m}_{1110} = \hat{m}_{1111} = 0,$$

since the last row sum is a sufficient statistic. We are hence left with a system of equations in twelve indeterminates which we call the *simplified four-cycle model*. This system consists of five binomials and eight linear equations. The following binomials are the minimal generators of the toric ideal of the simplified four-cycle model:

$$\hat{m}_{0011}\hat{m}_{1001} - \hat{m}_{0001}\hat{m}_{1011}$$
$$= \hat{m}_{0011}\hat{m}_{0110} - \hat{m}_{0010}\hat{m}_{0111}$$
$$= \hat{m}_{0001}\hat{m}_{0100} - \hat{m}_{0000}\hat{m}_{0101} = \hat{m}_{0010}\hat{m}_{1000} - \hat{m}_{0000}\hat{m}_{1010}$$
$$= \hat{m}_{0100}\hat{m}_{0111}\hat{m}_{1001}\hat{m}_{1010} - \hat{m}_{0101}\hat{m}_{0110}\hat{m}_{1000}\hat{m}_{1011} = 0.$$

We compute the marginal counts for the cliques in our model using the data set (4.13) and set these counts equal to the sum of the MLE for the cells associated with each clique as follows:

$$\hat{m}_{00++} = \hat{m}_{0000} + \hat{m}_{0001} + \hat{m}_{0010} + \hat{m}_{0011} = 4,$$
$$\hat{m}_{01++} = \hat{m}_{0100} + \hat{m}_{0101} + \hat{m}_{0110} + \hat{m}_{0111} = 4,$$
$$\hat{m}_{10++} = \hat{m}_{1000} + \hat{m}_{1001} + \hat{m}_{1010} + \hat{m}_{1011} = 5,$$
$$\hat{m}_{+00+} = \hat{m}_{0000} + \hat{m}_{0001} + \hat{m}_{1000} + \hat{m}_{1001} = 4,$$
$$\hat{m}_{++10} = \hat{m}_{0010} + \hat{m}_{0110} + \hat{m}_{1010} = 3,$$
$$\hat{m}_{++11} = \hat{m}_{0011} + \hat{m}_{0111} + \hat{m}_{1011} = 4,$$
$$\hat{m}_{0++1} = \hat{m}_{0001} + \hat{m}_{0101} + \hat{m}_{0011} + \hat{m}_{0111} = 4,$$
$$\hat{m}_{1++0} = \hat{m}_{1000} + \hat{m}_{1010} = 2.$$

Note that the following linear equations are implied and hence redundant in our system:

$$\hat{m}_{+01+} = \hat{m}_{0010} + \hat{m}_{0011} + \hat{m}_{1010} + \hat{m}_{1011} = 5,$$
$$\hat{m}_{+10+} = \hat{m}_{0100} + \hat{m}_{0101} = 2,$$
$$\hat{m}_{+11+} = \hat{m}_{0110} + \hat{m}_{0111} = 2,$$
$$\hat{m}_{++00} = \hat{m}_{0000} + \hat{m}_{0100} + \hat{m}_{1000} = 3,$$
$$\hat{m}_{++01} = \hat{m}_{0001} + \hat{m}_{0101} + \hat{m}_{1001} = 3,$$
$$\hat{m}_{0++0} = \hat{m}_{0000} + \hat{m}_{0100} + \hat{m}_{0010} + \hat{m}_{0110} = 4,$$
$$\hat{m}_{1++1} = \hat{m}_{1001} + \hat{m}_{1011} = 3.$$

The positive solution to these equations can be found numerically using *iterative proportional scaling*:

$$\begin{pmatrix} \hat{m}_{0000} & \hat{m}_{0001} & \hat{m}_{0010} & \hat{m}_{0011} \\ \hat{m}_{0100} & \hat{m}_{0101} & \hat{m}_{0110} & \hat{m}_{0111} \\ \hat{m}_{1000} & \hat{m}_{1001} & \hat{m}_{1010} & \hat{m}_{1011} \end{pmatrix} = \begin{pmatrix} 0.96 & 0.83 & 1.03 & 1.18 \\ 1.07 & 0.93 & 0.93 & 1.07 \\ 0.97 & 1.24 & 1.03 & 1.76 \end{pmatrix}.$$

Our main point, however, is to analyze the equations using symbolic algebra instead of numerical computation. We enter our five binomials and eight linear equations into the computer algebra system `Macaulay 2` by Grayson and Stillman [16]. To keep the notation simple, we replace the 12 indeterminates $\hat{m}_{0000}, \hat{m}_{0001}, \ldots, \hat{m}_{1011}$ by `a,b,...,l`. The command `gb MLE` computes the reduced Gröbner basis of our maximum likelihood equations in lexicographic term order:

```
i1 : R = QQ[a,b,c,d,e,f,g,h,i,j,k,l,
            MonomialOrder => Lex];

i2 : MLE = ideal(c*h-d*g, b*l-d*j, a*k-c*i,
                 a*f-b*e, e*h*j*k-f*g*i*l, a+b+c+d-4,
                 e+f+g+h-4, i+j+k+l-5, a+b+i+j-4,
                 c+g+k-3, d+h+l-4, b+f+d+h-4, i+k-2);

o2 : Ideal of R

i3 : gb MLE
```

The Gröbner basis consists of twelve polynomials:

$$\ell^5 - \frac{362}{39}\ell^4 + \frac{6713}{351}\ell^3 + \frac{110}{9}\ell^2 - \frac{2368}{39}\ell + \frac{480}{13},$$

$$k + \frac{6539}{22304}\ell^4 - \frac{58985}{33456}\ell^3 - \frac{513737}{602208}\ell^2 + \frac{490447}{100368}\ell - \frac{585}{2788},$$

$$j + \ell - 3, i + j + k + \ell - 5, h + \frac{1}{8}k^2 - \frac{1}{16}k\ell + \frac{1}{4}k - \frac{3}{16}\ell^2 + \frac{1}{16}\ell - \frac{7}{8},$$

$$g + h - 2, f - 2h - 2k + 3\ell - 2, e + f + g + h - 4, d + h + \ell - 4,$$

$$c + g + k - 3, b + f - \ell, a - f - g - h - k + 3.$$

The polynomials are in triangularized form; that is, each indeterminate is expressed in terms of indeterminates which come later in the alphabet. The only exception is the first equation in the Gröbner basis, which is a polynomial in the single variable $\ell$ and which we denote by $\psi(\ell)$. The properties of this polynomial which are relevant for our discussion are given by the following proposition:

PROPOSITION 3. *The polynomial $\psi(\ell)$ is irreducible over the rational numbers, and its Galois group is the symmetric group on five letters.*

Proposition 3 can be established using the computer algebra system `Maple` with the `Galois` command. The implications of this proposition are twofold. First, one cannot always find a rational solution to maximum likelihood for a nondecomposable model. Second, none of the five real roots of the equation can be expressed in terms of radicals. Thus, unlike a decomposable model, the MLE for a nondecomposable model cannot, in general, be given by an algebraic expression of the cell counts. For a more detailed discussion of Galois groups and irreducibility the reader can consult a book on Galois theory such as Stewart [27].

For the full four-cycle model, when none of the table counts is zero, the degree of the first equation in the Gröbner basis is 13 instead of 5. In particular, the polynomial of $\hat{m}_{1111}$ is irreducible of degree 13, and the other 15 coordinates $\hat{m}_{ijkl}$ of the maximum likelihood estimator are expressed as a polynomial with rational coefficients in $\hat{m}_{1111}$.

It would be interesting to find a combinatorial formula for the degree of the maximum likelihood estimator as a function of the structure of the undirected graphical model $A(G)$. A better understanding of this algebraic degree is likely to have applications in computational statistics.

4.7. *A characterization of decomposable models.* The following theorem provides a characterization of decomposable models.

THEOREM 4.4. *Let $G$ be an undirected graphical model for discrete variables. Then the following five statements are equivalent:*

(a) *$G$ is a decomposable graphical model.*

(b) *A distribution $P$ factors according to $G$ if and only if $P$ satisfies a set of quadratic binomials corresponding to global separation statements in $G$.*

(c) *The ideal $I_{A(G)}$ has a quadratic Gröbner basis in which each polynomial corresponds to a global separation statement in $G$.*

(d) *The maximum likelihood estimate for $G$ is a rational function.*

(e) *The set* image($\phi_{A(G)}$) *is closed.*

The fact that (a) implies (c) is Theorem 4.3. As described in Sections 4.3 and 4.6, it is known that (a) implies (d) and (a) implies (e). Note that (c) implies (b) is a trivial implication, so the only thing to prove for the above theorem is (b) implies (a), (d) implies (a) and (e) implies (a). We use constructions based on examples from previous sections to prove the result.

The essential idea of the proof is that every nondecomposable model contains a four-cycle and we prove all these claims by lifting the examples developed in Sections 4.5 and 4.6 for the four-cycle model to other nondecomposable models. For these results we use the following graph-theoretic definitions. A set of vertices $\mathbf{A}$ is *connected* in undirected graph $G$ if and only if there is a path in $G$ between every pair of vertices in $\mathbf{A}$ that only passes through vertices in $\mathbf{A}$. The sets $\mathbf{A}, \mathbf{B}, \mathbf{C}, \mathbf{D}, \mathbf{E}$ are a *nondecomposable partition for graph $G$ with vertices $\mathbf{X}$* if and only if (1) the sets $\mathbf{A}, \mathbf{B}, \mathbf{C}, \mathbf{D}, \mathbf{E}$ are disjoint, (2) $\mathbf{X} = \mathbf{A} \cup \mathbf{B} \cup \mathbf{C} \cup \mathbf{D} \cup \mathbf{E}$, (3) $\mathbf{A}, \mathbf{B}, \mathbf{C}, \mathbf{D}$ are not empty, (4) the subgraph of $G$ over $\mathbf{A}, \mathbf{B}, \mathbf{C}, \mathbf{D}$ is a cycle with no chords and (5) each of the sets $\mathbf{A}, \mathbf{B}, \mathbf{C}, \mathbf{D}$ is connected in $G$.

PROPOSITION 4. *If the undirected graph $G$ is nondecomposable, then there exists a nondecomposable partition for the graph.*

One can construct a nondecomposable partition for a nondecomposable graph $G$ with vertex set $\mathbf{X}$ as follows. First let $C_1, \ldots, C_n$ be a cycle in $G$ with no chords and length $n \geq 4$. One nondecomposable partition is given by the following five sets: $\mathbf{A} = \{C_1, \ldots, C_{i-1}\}$, $\mathbf{B} = \{C_i, \ldots, C_{j-1}\}$, $\mathbf{C} = \{C_j, \ldots, C_{k-1}\}$, $\mathbf{D} = \{C_k, \ldots, C_n\}$ where $1 < i < j < k \leq n$, and $\mathbf{E} = \mathbf{X} \setminus \{C_1, \ldots, C_n\}$. We are now prepared to prove the needed claims.

*Proof of* "(b) *implies* (a)." We use Example 8 to show that the zero set of $I_{A(G)}$ for an arbitrary nondecomposable model cannot be specified by quadratic binomials. We exhibit a probability distribution $p$ which is in the zero set of all quadratic binomials in $I_{A(G)}$ but is not in the zero set of $I_{A(G)}$.

Let $G$ be a nondecomposable graph and let the sets $\mathbf{A}, \mathbf{B}, \mathbf{C}, \mathbf{D}, \mathbf{E}$ form a nondecomposable partition of $G$. We define a probability distribution $p$ such that $p_{01001} = p_{01111} = p_{10011} = p_{10101} = 1/4$ where $p_{ijklm}$ is the probability that each variable in $\mathbf{A}$ has value $i$ and each variable in $\mathbf{B}$ has value $j$ and each variable in $\mathbf{C}$ has value $k$ and each variable in $\mathbf{D}$ has value $k$ and each variable in $\mathbf{E}$ has value $m$.

Consider the nonquadratic binomial $p_{01001}p_{01111}p_{10011}p_{10101} - p_{01011}p_{01101} \times p_{10001}p_{10111}$. This binomial lies in $I_{A(G)}$ because the intersection of any clique of $G$ and the set of vertices on the cycle defining the nondecomposable partition (i.e., $\mathbf{A} \cup \mathbf{B} \cup \mathbf{C} \cup \mathbf{D}$) is either the empty set, a singleton, or pair of adjacent vertices. Restricting the indices of the two quartic monomials to any clique gives two identical monomials, which means that the binomial lies in $I_{A(G)}$.

We claim that every quadratic binomial in $I_{A(G)}$ vanishes at $p$. Suppose not. Then there exists a binomial $p_a p_b - p_c p_d$ which lies in $I_{A(G)}$ and, after relabeling, our probability distribution $p$ satisfies $p_a = p_b = 1/4$ and $p_c p_d = 0$. Hence $a$ and $b$ are among the four basic events with positive probability. For any such pair $a, b$, it is easy to check that the sum of the two columns of $A(G)$ indexed by $a$ and $b$ cannot be written in any other way as a sum of columns of $A(G)$. The reason is that $a$ and $b$ agree in a connected subset of the $k$-cycle and they also disagree in a connected subset of the $k$-cycle. We conclude that every quadratic binomial in $I_{A(G)}$ vanishes at our probability distribution $P$, which completes the proof of the implication from (b) to (a).

*Proof of* "(d) *implies* (a)."    We lift the example described by (4.13) that demonstrates the potential nonrationality of the estimates for a four-cycle (Section 4.6) to an arbitrary nondecomposable graph. Let $G$ be an arbitrary nondecomposable graph and let the sets $\mathbf{A}, \mathbf{B}, \mathbf{C}, \mathbf{D}, \mathbf{E}$ form a nondecomposable partition of $G$. We define a data set for the variables in $G$ by expanding the data set defined by (4.13). Let $n_{ijklm}$ denote the count of cases in which all of the variables in $\mathbf{A}$ have the value $i$ and all of the variables in $\mathbf{B}$ have the value $j$ and so on. If we let the data set be

$$(4.14) \quad \begin{pmatrix} n_{00001} & n_{00011} & n_{00101} & n_{00111} \\ n_{01001} & n_{01011} & n_{01101} & n_{01111} \\ n_{10001} & n_{10011} & n_{10101} & n_{10111} \\ n_{11001} & n_{11011} & n_{11101} & n_{11111} \end{pmatrix} = \begin{pmatrix} 1 & 1 & 1 & 1 \\ 1 & 1 & 1 & 1 \\ 1 & 1 & 1 & 2 \\ 0 & 0 & 0 & 0 \end{pmatrix},$$

where all of the counts not shown are zero, then the estimate for $\hat{n}_{11111}$ will be identical to the nonrational maximum likelihood estimate of $\hat{m}_{1111}$ from Section 4.6.

*Proof of* "(e) *implies* (a)."    Let $G$ be a nondecomposable graph and let the sets $\mathbf{A}, \mathbf{B}, \mathbf{C}, \mathbf{D}, \mathbf{E}$ form a nondecomposable partition of $G$. We construct a sequence of distributions that factor according to $G$ but whose limiting distribution does not factor. We do so by defining pairwise potential functions for each of the edges in a graph $G$—that is, a log-linear distribution where every generator is a set of at most two variables. A log-linear distribution of this form will also factor according to $G$ because the pairwise potential functions can be combined into clique potentials. First $\psi_{\varnothing}(\cdot) = 1$. We consider pairs of vertices $\{X, Y\}$ connected by an edge in $G$. If $X \in \mathbf{E}$ or $Y \in \mathbf{E}$, then $\psi_{X,Y}(\cdot, \cdot) = 1$. If $\{X, Y\} \subset \mathbf{A}$, $\{X, Y\} \subset \mathbf{B}$, $\{X, Y\} \subset \mathbf{C}$ or $\{X, Y\} \subset \mathbf{D}$, then we define $\psi_{X,Y}(0,0) = \psi_{X,Y}(1,1) = n$ and otherwise $\psi_{X,Y}(\cdot, \cdot) = 1$. From the definition of a nondecomposable partition, $\mathbf{A}$ is connected to exactly two of the sets $\mathbf{B}, \mathbf{C}, \mathbf{D}$. Without loss of generality suppose that $\mathbf{A}$ is connected to $\mathbf{B}$ and $\mathbf{D}$. Finally we add the potentials for edges between the sets $\mathbf{A}, \mathbf{B}, \mathbf{C}, \mathbf{D}$. If $X \in \mathbf{A}$ and $X \in \mathbf{B}$, then $\psi_{X,Y}(x,y) = n^{(xy-y)}$ if $x \in \{0,1\}$ and $y \in \{0,1\}$ and $\psi_{X,Y}(x,y) = 1$ otherwise. If $X \in \mathbf{B}$ and $X \in \mathbf{C}$, then $\psi_{X,Y}(x,y) = n^{(xy-y)}$ if $x \in \{0,1\}$ and $y \in \{0,1\}$ and $\psi_{X,Y}(x,y) = 1$ otherwise. If $X \in \mathbf{C}$ and $X \in \mathbf{D}$, then $\psi_{X,Y}(x,y) = n^{(xy)}$ if $x \in \{0,1\}$ and $y \in \{0,1\}$ and $\psi_{X,Y}(x,y) = 1$ otherwise. If $X \in \mathbf{A}$ and $X \in \mathbf{D}$, then  $\psi_{X,Y}(x,y) = n^{(-xy)}$ if $x \in \{0,1\}$ and $y \in \{0,1\}$ and $\psi_{X,Y}(x,y) = 1$ otherwise. We consider the sequences of distributions defined by these pairwise potentials as $n \to \infty$. If we consider the four-cycle graph, then the limiting distribution is equal to the distribution given in Example 7 and the distribution does not have $A$-feasible support. In the limiting distribution for a general nondecomposable graph $G$, the variables in $\mathbf{E}$ are mutually independent and independent of all other variables and all of the variables within either $\mathbf{A}, \mathbf{B}, \mathbf{C}$ or $\mathbf{D}$ are deterministically related. Thus, checking whether the support of the limiting distribution is $A$-feasible reduces to the problem of checking the support of Example 7.

## APPENDIX: PROOFS OF THEOREMS 3.1 AND 3.2

We first note via the next lemma that it would be equivalent in the definition of the nonnegative toric variety $X_A$ to allow $u_1, \ldots, u_m, v_1, \ldots, v_m$ to be nonnegative real numbers rather than integers.

LEMMA A.1. *Let $Z_A$ be the set of all vectors $(x_1, \ldots, x_m) \in \mathbf{R}_{\geq 0}^m$ which satisfy*

$$(A.1) \qquad x_1^{u_1} x_2^{u_2} \cdots x_m^{u_m} = x_1^{v_1} x_2^{v_2} \cdots x_m^{v_m}$$

*whenever $u = (u_1, \ldots, u_m)$ and $v = (v_1, \ldots, v_m)$ are vectors of nonnegative real numbers which satisfy the $d$ linear relations*

$$(A.2) \qquad u_1 a_1 + u_2 a_2 + \cdots + u_m a_m = v_1 a_1 + v_2 a_2 + \cdots + v_m a_m.$$

*Then $Z_A = X_A$.*

PROOF. Clearly, $Z_A \subseteq X_A$. For the converse, let $x$ be a point in $X_A$. We need to show that (A.1) with $u, v$ being nonnegative real vectors holds for the point $x$.

A vector $b$ is *sign-compatible* with a vector $c$ if every nonzero entry of $b$ agrees in sign with the vector $c$. We denote by $c^+$ a vector whose $j$th entry equals $c_j$ for all nonnegative entries of $c$ and is zero otherwise. Similarly, we denote by $c^-$ a vector whose $j$th entry equals $-c_j$ for all negative entries of $c$ and is zero otherwise. Clearly, $c^+$ and $c^-$ are nonnegative vectors and $c = c^+ - c^-$.

From Lemma 4.10 of [28], there exist integer vectors $w_j$ that are sign-compatible with $w := u - v$ such that $w = \sum_j \alpha_j w_j$, where (i) $w_j^+$ and $w_j^-$ satisfy (A.2) and (ii) $\alpha_j \geq 0$. From (i) and the definition of $X_A$, we have $x^{w_j^+} = x^{w_j^-}$ for all $x \in X_A$. From (ii) and the fact that all of the $w_j$ are sign-compatible with $w$, we can write $w^+ = \sum_j \alpha_j w_j^+$ and $w^- = \sum_j \alpha_j w_j^-$. Because $0 \leq u - w^+ = v - w^-$, the expression $x^{u-w^+} = x^{v-w^-}$ is well defined and holds for all $x \in \mathbf{R}_{\geq 0}^m$. Therefore we can validly write $x^u = x^{w^+} x^{u-w^+}$ and $x^v = x^{w^-} x^{v-w^-}$. As $w^+ = \sum_j \alpha_j w_j^+$ and $w^- = \sum_j \alpha_j w_j^-$, and (A.1) holds for each pair $(w_j^+, w_j^-)$, it is straightforward to show that $x^{w^+} = x^{w^-}$, thus $x^u = x^v$. □

THEOREM 3.1. *A probability distribution $P$ factors according to $A$ if and only if $P$ lies in the nonnegative toric variety $X_A$ and the support of $P$ is $A$-feasible.*

PROOF. The only-if direction has been proved in Lemmas 1 and 2.

For the if direction, fix any vector $P \in X_A$ whose support $F = \operatorname{supp}(P)$ is $A$-feasible. We must prove that $P$ lies in the image of $\phi_A$, or equivalently that the system of (3.1) has a nonnegative real solution vector $(t_1, \ldots, t_d)$. Note that for the definition of $X_A$ we use nonnegative real exponents in (A.1) as justified by Lemma A.1.

Consider the following system of equations for the indeterminates $t_1, \ldots, t_d$:

$$(A.3) \qquad \prod_{i=1}^d t_i^{a_{ij}} = p_j > 0 \qquad \text{for } j \in F.$$

We claim that this system has a solution $(t_1, \ldots, t_d)$ all of whose coordinates are positive real numbers. Introducing new variables $\tau_i = \log(t_i)$, our claim is equivalent to the assertion that the following system of linear equations in $\tau_1, \ldots, \tau_d$ has a solution:

$$(A.4) \qquad \sum_{i=1}^d a_{ij}\tau_i = \log(p_j) \qquad \text{for } j \in F.$$

We proceed by contradiction.

Suppose that the system in (A.4) has no solution. This is a linear system of $|F|$ equations (over the field of the real numbers) in $d$ variables which can be written as $By = c$ where $B$ is an $|F| \times d$-matrix, and $c = (\log(p_j), j \in F)$ is a vector of length $|F|$. Assuming that (A.4) has no solution means that $c$ is not in the column space of $B$. Thus, there exists a row vector $q$ of length $|F|$ such that $qB$ is the zero vector but the inner product between the vectors $q$ and $c$ is not zero.

We now set $u_j = \max\{0, q_j\}$ and $v_j = \max\{0, -q_j\}$. Then $q_j = u_j - v_j$ and the identity $\sum_{j \in F} q_j\, a_{ij} = qB = 0$ translates into an identity of the form in (A.2) where $u_j = v_j = 0$ for all indices $j$ not in $F$. It follows from $\sum_{j \in F} q_j\, \log(p_j) \neq 0$ that

$$(A.5) \qquad \sum_{j \in F} u_j \log(p_j) \neq \sum_{j \in F} v_j \log(p_j).$$

Therefore,

$$\prod_{j \in F} p_j^{u_j} \neq \prod_{j \in F} p_j^{v_j}.$$

Consequently, the point $P$ does not satisfy (A.1) as required by all points on $X_A$. Hence, $P$ cannot lie on the nonnegative toric variety $X_A$, contrary to our assumption. Therefore the system in (A.4) can be solved for $\tau_i$, and hence (A.3) has a solution $t_i = \exp(\tau_i)$, $i \in \{1, \ldots, d\}$, in the positive reals.

The solution $(t_1, \ldots, t_d)$ just obtained is arbitrary at each index $i$ not in $I = \bigcup_{l \in F} \mathrm{supp}(a_l)$ because for each such $i$, $a_{ij} = 0$ for every $j \in F$. We now set $t_i = 0$ for all $i \in \{1, \ldots, d\} \setminus I$. Since $F$ is $A$-feasible, for each $j \notin F$ there exists an $i \in \{1, \ldots, d\} \setminus I$ such that $a_{ij} > 0$. Hence, $\prod_{i=1}^{d} t_i^{a_{ij}} = 0$ for $j \notin F$. We conclude that the modified vector $(t_1, \ldots, t_d)$ satisfies (3.1), and hence $P \in \mathrm{image}(\phi_A)$.   □

We fix the $d \times m$ matrix $A$ with columns $a_1, \ldots, a_m$ as before. A subset $F$ of $\{1, 2, \ldots, m\}$ is said to be *facial* if there exists a vector $c$ in $\mathbf{R}^d$ such that

(A.6)      $c^T a_i = 0$      for $i \in F$   and   $c^T a_i \geq 1$      for $i \in \{1, \ldots, m\} \setminus F$.

Hence, the vector $c$ is orthogonal to the columns whose index is in $F$ and not orthogonal to all other columns of $A$. The *characteristic vector* of $F$ is $(z_1, \ldots, z_m)$ with $z_i = 1$ if $i \in F$ and $z_i = 0$ if $i \notin F$.

LEMMA A.2.   *For a subset $F$ of $\{1, \ldots, m\}$ and a matrix $A$, the following statements are equivalent:*

(a) *$F$ is facial for $A$.*
(b) *The characteristic vector of $F$ lies in the nonnegative toric variety $X_A$.*
(c) *There exists a vector with support $F$ in the nonnegative toric variety $X_A$.*

PROOF.   Assume that (a) holds. We will first show that no nonzero nonnegative combination of the $a_i$, $i \notin F$, can be written as a linear combination of the $a_i$, $i \in F$. Let $c$ satisfy (A.6), and suppose $\sum_{i \notin F} \alpha_i a_i = \sum_{i \in F} \beta_i a_i$, where $\alpha_i \geq 0$, $i \notin F$. Then $0 \leq \sum_{i \notin F} \alpha_i \leq \sum_{i \notin F} \alpha_i c^T a_i = c^T \cdot (\sum_{i \notin F} \alpha_i a_i) = c^T \cdot (\sum_{i \in F} \beta_i a_i) = 0$, hence $\alpha_i = 0$ for $i \notin F$. Thus, there is no identity in (A.2) where $\mathrm{supp}(u) \subseteq F$ and $\mathrm{supp}(v)$ has nonempty intersection with $\{1, \ldots, m\} \setminus F$. Consequently, for every linear relation in (A.2), either both $\mathrm{supp}(u)$ and $\mathrm{supp}(v)$ are subsets of $F$ or neither of $\mathrm{supp}(u)$ and $\mathrm{supp}(v)$ is a subset of $F$. However, this means that the characteristic vector of the set $F$ satisfies (A.1) [namely, both sides of (A.1) are 0 or both sides are 1] whenever (A.2) holds. Equivalently, the characteristic vector of $F$ lies in $X_A$. Hence (a) implies (b).

Clearly, (b) implies (c). It remains to show that (c) implies (a). For this step we apply Farkas' lemma (linear programming duality); see Corollary 7.1e, Section 7.3, in [25]. Farkas' lemma reads as follows: *Let $D$ be a matrix* *and $e$ be a vector. Then the system $Dx \leq e$ has a solution $x$ if and only if $ye \geq 0$ for each nonnegative row vector $y$ with $yD = 0$.*

We define a matrix $D$ with $m + |F|$ rows as follows: The first $|F|$ rows are the vectors $-a_i$ for $i \in F$. The next $|F|$ rows are the vectors $+a_i$ for $i \in F$. The last $m - |F|$ rows are the vectors $-a_i$ for $i \notin F$. Let $e$ be the column vector with $m + |F|$ coordinates as follows: The first $|F|$ entries are 0. The next $|F|$ entries are 0. The last $m - |F|$ entries are $-1$.

Suppose $x_0 \in X_A = Z_A$, $\mathrm{supp}(x_0) = F$. Any nonnegative row vector $y = (y^{(1)}, y^{(2)}, y^{(3)})$ of respective lengths $(|F|, |F|, m - |F|)$ satisfying $yD = 0$ must have $y^{(3)} = 0$. Otherwise, taking $u = (y^{(2)}, 0)$ and $v = (y^{(1)}, y^{(3)})$ would satisfy (A.2) but contradict (A.1) for $x_0$, since $x_0^u > 0$ whereas $x_0^v = 0$. But then, when $y^{(3)} = 0$, each nonnegative solution of $yD = 0$ trivially satisfies $ye \geq 0$, hence by Farkas' lemma, $Dx \leq e$ has a solution. Consequently, (a) holds. □

THEOREM 3.2.  *A probability distribution $P$ factors according to $A$ or is the limit of probability distributions that factor according to $A$ if and only if $P$ lies in the nonnegative toric variety $X_A$.*

PROOF.  The claim is that $X_A = \mathrm{closure}(\mathrm{image}(\phi_A))$. By Lemma 2, the image of $\phi_A$ lies in $X_A$. The set $X_A$ is closed in $\mathbf{R}_{\geq 0}^m$ because it is defined by polynomial equations. Hence the closure of $\mathrm{image}(\phi_A)$ is contained in $X_A$. Therefore it suffices to prove that $X_A \subseteq \mathrm{closure}(\mathrm{image}(\phi_A))$. This is shown by taking a point $P \in X_A \setminus \mathrm{image}(\phi_A)$ and showing that $P$ lies in the closure of $\mathrm{image}(\phi_A)$. The argument is composed of three steps. Given a point $P \in X_A \setminus \mathrm{image}(\phi_A)$, we first define a sequence of points $P(\varepsilon)$; we then prove that $\lim_{\varepsilon \to 0} P(\varepsilon) = P$; and finally, we prove that $P(\varepsilon) \in \mathrm{image}(\phi_A)$ for all $\varepsilon > 0$.

Let $P \in X_A \setminus \mathrm{image}(\phi_A)$ and $F = \mathrm{supp}(P)$. In order to define $P(\varepsilon)$, consider the following system of equations for the indeterminates $t_1, \ldots, t_d$:

$$(A.7) \qquad \prod_{i=1}^{d} t_i^{a_{ij}} = p_j > 0 \qquad \text{for } j \in F.$$

This system of equations is identical to (A.3). We have shown that it has a solution $(t_1, \ldots, t_d)$ all of whose coordinates are positive real numbers. By Lemma A.2, the set $F = \mathrm{supp}(P)$ is facial. We now fix $c \in \mathbf{R}^d$ so that (A.6) holds. We introduce a positive real parameter $\varepsilon > 0$, and define the vector

$$P(\varepsilon) = \left( \varepsilon^{c^T a_1} \prod_{i=1}^{d} t_i^{a_{i1}}, \varepsilon^{c^T a_2} \prod_{i=1}^{d} t_i^{a_{i2}}, \ldots, \varepsilon^{c^T a_m} \prod_{i=1}^{d} t_i^{a_{im}} \right).$$

The condition in (A.6) implies that $\lim_{\varepsilon \to 0} P(\varepsilon) = P$ because for every $j$ in $F$, $\varepsilon^{c^T a_j}$ is always 1 and so $p_j$ tends to $\prod_{i=1}^{d} t_i^{a_{ij}}$, and for every $j$ not in $F$,

$\varepsilon^{c^T a_j}$ tends to 0 and so $p_j$ tends to 0. Finally, note

$$\exp(c^T a_j) = \exp(c_1 a_{1j} + c_2 a_{2j} + \cdots + c_d a_{dj})$$
$$= \exp(c_1)^{a_{1j}} \exp(c_2)^{a_{2j}} \cdots \exp(c_d)^{a_{dj}}$$
$$= \prod_{i=1}^{d} \exp(c_i)^{a_{ij}},$$

where exp denotes the exponential function with base $\varepsilon$. Thus, the $j$th coordinate of the vector $P(\varepsilon)$ equals $\prod_{i=1}^{d}(\exp(c_i)^{a_{ij}} t_i^{a_{ij}}) = \prod_{i=1}^{d}(\exp(c_i)t_i)^{a_{ij}}$. Hence $P(\varepsilon)$ is the image of the strictly positive vector $(\exp(c_1)t_1, \exp(c_2)t_2, \ldots, \exp(c_d)t_d)$ under the map $\phi_A$. This shows that $P(\varepsilon)$ lies in the image of $\phi_A$ for all $\varepsilon > 0$.  $\square$

**Acknowledgments.** We thank the anonymous referees and especially an Associate Editor whose comments on previous manuscripts have significantly improved the final version of this paper.

This work was done in part while the first-named author visited Microsoft Research.

D. GEIGER
COMPUTER SCIENCE DEPARTMENT
TECHNION
HAIFA 36000
ISRAEL
E-MAIL: dang@cs.technion.ac.il

C. MEEK
MICROSOFT RESEARCH
REDMOND, WASHINGTON 98052
USA
E-MAIL: meek@microsoft.com

B. STURMFELS
DEPARTMENT OF MATHEMATICS
UNIVERSITY OF CALIFORNIA
BERKELEY, CALIFORNIA 94720
USA
E-MAIL: bernd@math.berkeley.edu